\documentstyle[12pt]{article}
\font\tenmsb=msbm10
\newfam\msbfam
\textfont\msbfam=\tenmsb
\def\Bbb#1{{\fam\msbfam#1}}

\title{\bf A new proof of the Abhyankar-Moh-Suzuki Theorem}
\date{\ }
\begin{document}
\author {Leonid Makar-Limanov}

\newcommand{\C}{\Bbb C}
\newcommand{\Z}{\Bbb Z}
\newcommand{\de}{\partial}
\newcommand{\J}{\mbox{J}}
\newcommand{\mod}{\mbox{mod}}
\newcommand{\Def}{\mbox{def}}
\newcommand{\gap}{\mbox{gap}}
\newcommand{\spn}{\mbox{span}}
\newcommand{\chr}{\mbox{char}}

\maketitle
\vskip -1cm
\hskip 5cm{{\it To the memory of Shreeram Abhyankar  }}

\hskip 5cm{{\it  whose sudden death was a profound shock}}

\hskip 5cm{{\it and a tremendous loss}}
\vskip 0.5cm

\begin{abstract} This note contains a complete proof of the
Abhyankar-Moh-Suzuki theorem (in characteristic zero case).
\end{abstract}

\textbf{Introduction.}\\

In the zero characteristic case the AMS Theorem which was independently proved by
Abhyankar-Moh and Suzuki (see [AM] and [Su]) and later reproved by
many authors (see [AO], [AB], [Es], [Gu], [GM], [Ka], [Mi], [Ri], [Ru], [Za]; the list is probably incomplete) states the following

\textbf{AMS Theorem.}  If $f$ and $g$ are polynomials in $K[z]$
of degrees $n$ and $m$ for which $K[f, g] = K[z]$
then $n$ divides $m$ or $m$ divides $n$.\\

Here is the plan of a proof. We start with an algorithm which produces the monic irreducible dependence for any pair of polynomials $f, g \in K[z]$ where $K$ is a field of any
characteristic. This algorithm also produces a \textit{standard} linear basis of $K(f)[g]$ over $K(f)$ which consists of elements of $K[f, f^{-1}, g]$ of pairwise different degrees.  When characteristic is zero or when characteristic does not divide the degree of $g$ the standard basis consists of polynomials from $K[f,g]$ monic in $g$. After this is established the AMS Theorem follows almost immediately.\\

\textbf{Irreducible dependence of two polynomials.}\\

In this section we describe an algorithm for finding the minimal algebraic
dependence between $f, \ g \in K[z]$ where $K$ is a field of any
characteristic. The algorithm seems to be new though it is
not very different from the algorithm suggested by David Richman and
Barbara Peskin (see [PR], [R], [Es],  and [Ka]). In fact, when $m$ and $n$ are
relatively prime this is the algorithm from [PR] but when $m$ and $n$ are not
relatively prime the algorithm from [PR] requires more intermediate steps.

Let $E = K(z)$ and $F = K(f(z))$ be the fields of rational
functions in $z$ and $f(z)$ correspondingly. Since $F \subset E$
we can consider $E$ as a vector space over $F$. Denote by $[E:F]$
the dimension of this vector space.\\

The next two Lemmas may be skipped by a reader who knows that
there exists an irreducible polynomial dependence between $f$ and
$g$ which is given by a polynomial monic in $g$.\\

\textbf{Lemma 1.} $[E:F] = n = \deg(f)$ and $\{1, z, \dots,
z^{n-1}\}$ is a basis of $E$ over $F$.\\
\textbf{Proof.} The degrees of $\alpha_iz^i$
where $\alpha_i \in K[f(z)]$ and $0 \leq i < n$ are different for
different $i$'s. Hence the elements $\{1, z, \dots, z^{n-1}\}$ are
linearly independent over $F$. If $[E:F] > n$ take $n+1$ elements linearly
independent over $F$ and multiply them by a common denominator to obtain
$n+1$ elements of $K[z]$ linearly independent over $F$.
On the other hand $K[z] = \bigoplus\limits_{i =
0}^{n-1}z^iK[f(z)]$ since for any non-negative
$k$ a monomial $z^k$ is contained in $\bigoplus\limits_{i =
0}^{n-1}z^iK[f(z)]$.
Hence $K[z]$ cannot contain $n + 1$ elements linearly independent over $F$.$\Box$\\

Let $g \in K[z]$. By the previous Lemma there exists a non-trivial
relation $\sum\limits_{i = 0}^{n} \alpha_i g^i = 0$, i.e. there
exists a non-zero element $P(x, y) \in A = K(x)[y]$ for which
$P(f(z), g(z)) = 0$. We will assume that $k = \deg_y(P)$ is
minimal possible and that $P$ is monic in $y$. Then $P$ is an irreducible element of $A$ and if
$Q(f, g) = 0$ for some $Q \in A$ then $Q$ is divisible by $P$ (by the Euclidean algorithm).\\

\textbf{Lemma 2.} $P \in K[x, y]$.\\
\textbf{Proof.} Since $P = y^k + \sum\limits_{i = 0}^{k-1} p_i(x)
y^i$ where $p_i \in K(x)$ we can multiply $P$ by the least common
denominator $D(x) \in K[x]$ of $p_i$ and obtain a polynomial $DP
\in K[x, y]$ which is irreducible in $K[x, y]$. In order to prove
that $D = 1$ it is sufficient to find an element $Q \in K[x, y]$
such that $Q(f, g) = 0$ and $Q$ is monic in $y$. Indeed, $Q$ must
be divisible by $DP$ in $K[x, y]$ by the Gauss lemma, which is
possible only if $D = 1$.

For a natural number $l$ define $Q_l \in K[x, y]$ as
$Q_l = y^l + R_l$ where $\deg_y(R_l) < l$ and $\deg_z(Q_l(f,g))$ is the minimal possible.
Put $e_l = \deg_z( Q_l(f,g))$ when $Q_l(f,g) \neq 0$. If $a > b$ and $e_a \equiv e_b \pmod{n}$
then $e_a < e_b$ because otherwise we can find $j \in \Z^+$ and $c
\in K$ so that $\deg_z(Q_a(f,g) - c f^j Q_b(f,g)) <
\deg_z(Q_a(f,g))$. Therefore we can have only a finite number of
$e_a$ which means that $Q_a(f, g) = 0$ for a
sufficiently large $a$.$\Box$\\

Let us describe now a procedure which will produce $P$.
First an informal description. Raise $g$ to the smallest possible
power $a$ so that by subtracting some power of $f$ (with an
appropriate coefficient) the degree of $g^a$ can be decreased. If the
result has the degree which can be decreased by subtracting a
monomial in $f$ and $g$, do it and continue until the degree of
the result cannot be decreased. Since different monomials in $f$
and $g$ can have the same degree, use only monomials with power of
$g$ less than $a$. Then the choice of a monomial with given degree
is unique. If the result $h$ is zero it gives the dependence we
are looking for. If not, raise $h$ to the smallest possible power $a'$
so that the degree of $h^{a'}$ can be decreased by subtracting a monomial
in $f, \ g$ and on further steps use for reduction purposes the
monomials in $f, \ g, \ h$ with appropriately restricted
powers of $g$ and $h$. After several steps like that an
algebraic dependence will be obtained.

It is easy to implement
this procedure and it works nicely on examples. On the other hand
why should it stop? If a monomial with a negative power of $f$ is used at some stage, we obtain a rational function and it is not clear why the process stops after a finite number of the degree reductions.
Also even if all monomials which are used in reductions have $f$ in positive power, and then it is clear that every step stops after a finite number of reductions of the degree, since the degrees from a step to a step may increase, why the number of steps is finite?

Here is an example where negative powers of $f$ appear. Take $f = z^4$, $g = z^6 - z$. We have to
start with $g^2 - f^3 = -2z^7 + z^2$ and $h = -2z^7 + z^2$. Next,
$h^2 - 4f^2g = z^4$ and $h^2 - 4f^2g - f = 0$. So $(g^2 - f^3)^2 -
4f^2g - f = 0$.

Assume now that the ground
field has characteristic $2$. Then $g^2 - f^3 = z^2$ and we can proceed with
the degree reduction to get $h = g^2 - f^3 - f^{-1}g = z^{-3}$ and a dependence $h^2 - f^{-3}g -
f^{-2}h = 0$
in which miraculously all negative
powers of $f$ disappear: $h^2 - f^{-3}g -
f^{-2}h = g^4 - f^6 -
f^{-2}g^2 - f^{-3}g - f^{-2}g^2 - f - f^{-3}g = g^4 - f^6 - f$.\\

\textsc{formal description.}\\

Below $\deg(h)$ denotes the $z$-degree of $h \in K(z)$ defined as the difference of the degrees of the numerator and the denominator of $h$.\\

\emph{First step.}\\

Put $g_0 = g$. Let $\deg(g_0) = m_0$ and $\deg(f) = n$. Find the
greatest common divisor $d_0$ of $n$ and $m_0$. Take
the smallest positive integers $a_0 = {n\over d_0}, \ b_0 = {m_0\over d_0}$ for which $\deg(g_0^{a_0}) = \deg(f^{b_0})$. Find $k_0 \in K$ for which
$m_{0,1} = \deg(g_0^{a_0} - k_0f^{b_0}) < \deg(g_0^{a_0})$. If
$m_{0,1}$ is divisible by $d_0$ find a monomial $f^ig_0^{j_0}$
with $0 \leq j_0 < a_0$ and $\deg(f^ig_0^{j_0}) = m_{0,1}$, find
$k_1 \in K$ for which $m_{0,2} = \deg(g_0^{a_0} - k_0f^{b_0} -
k_1f^ig_0^{j_0}) < m_{0,1}$ and so on.\\

If the procedure does not stop we failed.

If after a finite number of reductions $m_{0, i}$ which is not
divisible by $d_0$ is obtained, denote the corresponding
expression by $g_1$ and make the next step.

If after a finite number of reductions zero is obtained, we have a
dependence and stop. \\

\emph{ Generic step.}\\

Assume that after $s$ steps we obtained $g_0, \dots, g_s$ where
$g_s \neq 0$. Denote $\deg(g_i)$ by $m_i$ and the greatest common divisor  $(n, m_0,
\dots, m_i)$ of $n, m_0, \dots, m_i$ by $d_i$. The numbers $d_i$ are positive while $m_i$
can be negative. Put $d_{-1} = n$ and $a_i = \frac{d_{i-1}}{d_i}$
for $0 \leq i \leq s$. (Clearly $a_sm_s$ is divisible by $d_{s-1}$
and $a_s$ is the smallest integer with this property.) Call a
monomial ${\bf m} = f^ig_0^{j_0} \dots g_s^{j_s}$ with $0 \leq j_k
< a_k$ $s$-\textit{standard}.

Find an $s-1$-standard monomial ${\bf m}_{s,0}$ with $\deg({\bf
m}_{s,0}) = a_sm_s$ and $k_0 \in K$ for which $m_{s,1} =
\deg(g_s^{a_s} - k_0 {\bf m}_{s,0}) < a_sm_s$.
If $m_{s,1}$ is divisible by $d_s$ find an $s$-standard monomial
${\bf m}_{s,1}$ with $\deg({\bf m}_{s,1}) = m_{s,1}$ and $k_1 \in
K$ for which $m_{s,2} = \deg(g_s^{a_s} - k_0 {\bf m}_{s,0} - k_1
{\bf m}_{s,1})
< m_{s,1}$ and so on.
(We will check in Lemma 3 that any number divisible by $d_s$ is the degree of an $s$-standard monomial.)\\

If the procedure does not stop we failed.

If after a finite number of reductions $m_{s, i}$ which is not
divisible by $d_s$ is obtained, denote the corresponding
expression by $g_{s+1}$ and make the next step.

If after a finite number of reductions zero is obtained, we have a
dependence and stop. \\

\textbf{Remark.} If $g_{i+1}$ is constructed then $d_{i+1} = (d_i, m_{i+1}) < d_i$ since $m_{i+1}$ is not divisible by $d_i$; therefore $d_0 > d_1 > \dots, > d_s$. $\Box$\\

To prove that failure is not an option we should know more about $s$-standard monomials. In the sequel $g_i$ are considered as the elements of $L = K[f, f^{-1}, g]$ where $f, \ g$ are variables, as well as the elements of $E = K(z)$. \\

\textbf{Lemma 3.} If the elements $g_0, g_1, \dots, g_s$ are defined then\\
(a) Any number divisible by $d_s = (n,m_0, \dots, m_s)$ is the degree of an $s$-standard monomial and this monomial is uniquely defined;\\
(b) For any $d < a_s\deg_g(g_s)$ there exists an $s$-standard monomial $\bf m$ with $\deg_g({\bf m}) = d$.\\
\textbf{Proof.} In this proof $s$-standard monomials do not contain $f$.\\
(a) The degrees of different $s$-standard
monomials are different $\mod\, n$. Indeed, if $\sum\limits_{k =
0}^{s} j_k m_k \equiv \sum\limits_{k = 0}^{s} i_k m_k \, (\mod\,
n)$ then $j_s m_s \equiv i_s m_s \, (\mod\, d_{s-1})$. Therefore $j_s =
i_s$ because $0 \leq i_s, \ j_s < a_s$ and $|j_s -
i_s|m_s$ is not divisible by $d_{s-1}$ if $0 < |j_s - i_s| < a_s$ by the definition of $a_s$.
So $j_s = i_s$ and we can omit them from the sums and proceed to
prove that $j_{s-1} = i_{s-1}$, etcetera. There is $\prod\limits_{k =
0}^{s} a_k = {d_{-1}\over d_s} = {n\over d_s}$ different $s$-standard monomials and there is ${n\over d_s}$ residues $\mod\, n$ divisible by $d_s$. Hence any number divisible by $d_s$ is the degree of a unique $s$-standard monomial $f^i{\bf m}$.\\
(b) The elements $g_i \in L = K[f, f^{-1}, g]$. It is easy to
check by induction that $\deg_g(g_t) = a_0 \dots a_{t-1}$ for $t \leq s$. The
base $\deg_g(g_0) = 1$ is clear since $g_0 = g$. Assume that
$\deg_g(g_k) = a_0 \dots a_{k-1}$ for $k < t + 1$.
For a $t$-standard monomial ${\bf m} = g_0^{j_0} \dots g_t^{j_t}$ the degree
$\deg_g({\bf m}) = \sum\limits_{l = 0}^{t} j_l\deg_g(g_l) \leq
\sum\limits_{l = 0}^{t} (a_l - 1) \deg_g(g_l) = \sum\limits_{l = 0}^{t-1} (\deg_g(g_{l+1}) -  \deg_g(g_l)) + (a_t - 1) \deg_g(g_t) = \deg_g(g_t) - 1 + (a_t - 1) \deg_g(g_t) = a_t \deg_g(g_t) - 1$ under the induction assumption. Therefore $\deg_g({\bf m}) \leq a_t\deg_g(g_t) - 1$.
Now, $g_{t+1} = g_t^{a_t} - r_t(f, g_0, \dots, g_t)$.
Since all monomials of $r_t$ are $t$-standard, $\deg_g(r_t) \leq a_t\deg_g(g_t) - 1$ and $\deg_g(g_{t+1}) = \deg_g(g_t^{a_t}) = a_0 \dots a_{t-1}a_t$.

If ${\bf m} = g_0^{j_0} \dots g_s^{j_s}$ and $\deg_g({\bf m}) = \sum\limits_{k = 0}^{s} j_k \deg_g(g_k)  =
\sum\limits_{k = 0}^{s} i_k \deg_g(g_k)$ then $j_0 \equiv i_0 \,
(\mod\, a_0)$ and $j_0 = i_0$ because $0 \leq j_0 < a_0$ and $0
\leq i_0 < a_0$; we can proceed to prove that $j_1 = i_1$ since then $j_1 \equiv i_1 \,
(\mod\, a_1)$ etc.. Hence different $s$-standard monomials have different
$g$-degrees.
There is exactly $a_0 \dots a_s = a_s\deg_g(g_s)$ $s$-standard
monomials and $\deg_g({\bf m}) < a_s\deg_g(g_s)$ for $s$-standard monomials. Therefore
we have an $s$-standard monomial with $g$-degree equal to $d$ for any $d < a_s\deg_g(g_s)$.$\Box$\\

\textbf{Remark.} A standard monomial ${\bf m} = f^ig_0^{j_0} \dots
g_s^{j_s}$ is completely determined by $i$ and $\deg_g({\bf m})$. $\Box$\\

\textbf{Lemma 4.} If the elements $g_0, g_1, \dots, g_s \in K(z)$ are defined and $g_s \neq 0$ then $g_{s+1}$ is also defined. \\
\textbf{Proof.}
The field $E = K(z)$ is a vector space over its subfield $F = K(f(z))$.
Denote by $V_s$ the subspace of $E$ generated over $F$ by all $s$-standard monomials.
There are two possibilities: $g_s^{a_s} \not\in V_s$ and $g_s^{a_s} \in V_s$.

Since the degrees of different $s$-standard monomials not containing $f$ are
different $\mod\, n$ (see the proof of Lemma 3 (a)) they are linearly independent over $F$
and form a \textit{standard} basis $B_s$ of $V_s$.

Assume that $g_s^{a_s} \not\in V_s$. As we know $E$ is $n$-dimensional
over $F$ and $\{1, z, \dots, z^{n-1}\}$ is a basis of $E$ over $F$
(Lemma 1). The standard basis $B_s$ of $V_s$ contains $\prod\limits_{i=0}^s a_i =
{d_{-1}\over d_s} = {n\over d_s}$ elements. The degrees of the
elements of $B_s$ are divisible by $d_s$. The
elements
$\{z^i{\bf m_j} \, | \, 0 \leq i < d_s\}, \ {\bf m_j} \in B_s$ are linearly independent over $F$
since their degrees are different $\mod\, n$. Since there is $n$ of them
they form a basis of $E$ over $F$.
Write $g_s^{a_s} = \sum\limits_{{\bf m_j} \in B_s} \delta_{\bf j}{\bf m_j} + \sum\limits_{{\bf m_j} \in B_s}\sum\limits_{k = 1}^{d_s -1}
\epsilon_{k, \bf j} z^k{\bf m_j}$ where
$\delta_{\bf j}, \, \epsilon_{k, \bf j} \in F$.
The second sum is not zero and $D = \deg(\sum\limits_{{\bf m_j} \in B_s}\sum\limits_{k = 1}^{d_s -1}
\epsilon_{k, \bf j} z^k{\bf m_j})$ is not divisible
by $d_s$.

A rational function $\delta_{\bf j}$ can be approximated by a Laurent polynomial and written as
$\delta_{\bf j} = \sum\limits_{i = -N}^{M} c_{{\bf j},i} f^i +
r_{{\bf j},N}$ where $c_{{\bf j},i} \in K, \ r_{{\bf j},N} \in F, \ \deg(c_{{\bf j},i} f^i{\bf m_j}) > D$, and $\deg(r_{{\bf j},N}{\bf m_j}) < D$.
Therefore $g_s^{a_s} - \sum\limits_{{\bf m_j} \in B_s} \delta_{\bf j}{\bf m_j} = g_s^{a_s} - \sum\limits_{{\bf m_j} \in B_s}(\sum\limits_{i = -N}^{M} c_{{\bf j},i} f^i + r_{{\bf j},N}){\bf m_j}$ and $g_s^{a_s} - \sum\limits_{{\bf m_j} \in B_s}\sum\limits_{i = -N}^{M} c_{{\bf j},i} f^i{\bf m_j} =
\sum\limits_{{\bf m_j} \in B_s}(\sum\limits_{k = 1}^{d_s -1}\epsilon_{k, \bf j} z^k + r_{{\bf j},N}){\bf m_j}$ where $\deg(\sum\limits_{{\bf m_j} \in B_s}(\sum\limits_{k = 1}^{d_s -1}\epsilon_{k, \bf j} z^k + r_{{\bf j},N}){\bf m_j}) = \deg(\sum\limits_{{\bf m_j} \in B_s}\sum\limits_{k = 1}^{d_s -1}\epsilon_{k, \bf j} z^k{\bf m_j})$ is not divisible by $d_s$.
Hence $g_s^{a_s} - \sum\limits_{{\bf m_j} \in B_s}\sum\limits_{i = -N}^{M} c_{{\bf j},i} f^i{\bf m_j} = g_{s+1}$.

\smallskip
\smallskip

If $g_s^{a_s} \in V_s$ then $g_s^{a_s} = \sum\limits_{{\bf m_j} \in B_s} \delta_{\bf j}{\bf m_j}$ for some
$\delta_{\bf j} \in F$.
Let us show that in this case $g_{s+1} = 0$.
Recall that every $s$-standard monomial belongs to $L = K[f, f^{-1}, g]$.
Consider $P = g_s^{a_s} - \sum\limits_{{\bf m_j} \in B_s} \delta_{\bf j}{\bf m_j}$ as an element of $F[g]$.
By the proof of Lemma 3 (b) $\deg_g({\bf m_j}) < a_s\deg_g(g_s)$. Hence
$\deg_g(P) = a_s\deg_g (g_s)$ and $P$ is a monic polynomial in $g$.
Similarly, $g_i$ for $i \leq s$ and elements of $B_s$ are monic polynomials in $F[g]$. In Lemma 3 (b) we checked that $g$-degrees of elements of $B_s$ are pairwise different and that for any $d < a_s\deg_g(g_s)$ there is an element $b_d \in B_s$ with $\deg_g(b_d) = d$.
If $P$ is reducible in $F[g]$ then $P =
Q_1 Q_2$ where $\deg_g(Q_i) < \deg_g(P)$ and $Q_1, \  Q_2$ are non-zero elements of $F[g]$.
Hence $Q_1, \ Q_2$ can be presented as non-zero linear combinations (over $F$) of elements from $B_s$.
But $B_s$ is a basis of $V_s$ and $Q_i(f(z), g(z)) \neq 0$ while $P(f(z), g(z)) = 0$, a contradiction. Hence $P$ is irreducible and $P(f, g) \in K[f, g]$ by Lemma 2.
Now, $g_s^{a_s} \in L$ since $g_s \in L$.
Therefore $\sum\limits_{{\bf m_j} \in B_s} \delta_{\bf j}{\bf m_j} = g_s^{a_s} - P \in L$ and all $\delta_{\bf j} \in K[f, f^{-1}]$. (A presentation of $\sum\limits_{{\bf m_j} \in B_s} \delta_{\bf j}{\bf m_j}$ through the standard basis is unique since the elements of the standard basis have different $g$-degrees, also elements of $B_s$ are monic polynomials in $L$.) Consequently $\sum\limits_{{\bf m_j} \in B_s} \delta_{\bf j}{\bf m_j}$ can be presented as a finite sum of $s$-standard monomials with the coefficients from $K$ and the algorithm will produce zero after a finite number of steps. The monic irreducible relation $P(f, g)$ is also produced.
$\Box$\\

\textbf{Lemma 5.} After a finite number of steps the algorithm produces zero and a relation.\\
\textbf{Proof.} If the elements $g_0, \dots, g_{n+1}$ are defined and $g_{n+1} \neq 0$ then $\dim (V_{n+1}) > n$
since by the previous Lemma $\dim(V_i) < \dim (V_{i+1})$ if $g_{i+1} \neq 0$. But $\dim(V_i) \leq \dim(E) = n$. Hence $g_{s+1} = 0$ for some $s < n$ and $P = g_s^{a_s} - \sum\limits_{{\bf m_j} \in B_s} \delta_{\bf j}{\bf m_j}$
is a relation. $\Box$\\

So the algorithm works and we even know that $P \in L$ does not contain negative powers of $f$.\\

\textbf{Proof of AMS.}\\

Now we are ready to prove the AMS Theorem.

If $g_{s+1} = 0$ then by Lemma 3 (b) and since ${\bf m_j} \in B_s \subset K[f, f^{-1}, g]$ are elements monic in $g$, any element $h \in K[f, g]$ can be presented as a sum $h = \sum\limits_{{\bf m_j} \in B_s} \delta_{\bf j}{\bf m_j}$ where $\delta_{\bf j}(f) \in K(f)$.\\

\textbf{Lemma 6.} If characteristic of $K$ is zero then all $g_i$ are polynomials in $f$ and $g$.\\
\textbf{Proof.} Order the monomials $f^ig^j$ of $L = K[f, f^{-1},
g]$ lexicographically by $\deg_g, \deg_f$. Call a monomial
\textit{negative} if its $f$-degree is negative, otherwise call it
\textit{positive}. For an element $h \in L$ introduce a function
\textit{gap} as follows.
If $h \not\in K[f, g]$ then $\gap(h) = \overline{h} \div \widetilde{h}$ where
$\overline{h}$ is the largest monomial of $h$
and $\widetilde{h}$ is the largest negative monomial of $h$; if $h \in K[f, g]$ then $\gap(h) =
\infty$. Define $\infty$ to be larger than any monomial.

We will use the following properties of gap which are easy to
check:

(a) $\gap(h_1h_2) \geq \min(\gap(h_1), \gap(h_2))$;

(b) $\gap(h^d) = \gap(h)$ if $h$ is monic in $g$ and the
characteristic is zero;

(c) $\gap(fh) \geq \gap(h)$.\\

The plan is to show that $\gap(g_{j+1}) \leq \gap(g_j)$. Since we know that
the last $g_{s+1}$ which gives an irreducible dependence of $f(z)$ and
$g(z)$ is a polynomial in $f$ and $g$, this will imply that $\gap(g_j) = \infty$
for all $j$ and hence the Lemma because $\gap(h) = \infty$ is equivalent to $h \in K[f, g]$.

Let us use induction. The base of induction $\gap(g_1) \leq
\gap(g_0)$ is obvious since $\gap(g_0) = \infty$. Assume that
$\gap(g_{j+1}) \leq \gap(g_j)$ if $j < k$. If $g_k \in K[f, g]$
then $\gap(g_{k+1}) \leq \gap(g_k)$. So let $g_k \in L \setminus K[f, g]$

Since $g_{k+1}= g_k^{a_k} - r_k$ and $\gap(g_k^{a_k}) = \gap(a_k)$ it is sufficient to
check that the largest negative monomial of $r_k$ cannot cancel out
the largest negative monomial of $g_k^{a_k}$: then the largest negative monomial of $g_{k+1}$ is not smaller than the largest negative monomial of $g_k^{a_k}$ while their largest monomials are the same.

As above, call a $k$-standard monomial \textit{negative} if its
$f$-degree is negative and \textit{positive} otherwise. Let ${\bf
m} = f^ig_0^{j_0} \dots g_k^{j_k}$ be a $k$-standard monomial.
From the properties of gap mentioned above it follows that
$\gap(g_0^{j_0} \dots g_k^{j_k}) \geq \gap(g_k)$. Indeed
$\gap(g_i^{j_i}) = \gap(g_i)$ since $g_i$ is monic in $g, \ \gap(h_1h_2) \geq \min(\gap(h_1), \gap(h_2))$, and $\gap(g_i) \geq \gap(g_k)$ by the induction assumption.
Also if
$i \geq 0$ then $\gap(f^i h) \geq \gap(h)$, so $\gap({\bf m}) \geq
\gap(g_k)$ for a positive $k$-standard monomial ${\bf m}$. If $i <
0$ then $\gap({\bf m}) = 1$ since $g_0^{j_0} \dots g_k^{j_k}$ is
monic in $g$ and the largest monomial of ${\bf m} = f^ig_0^{j_0} \dots g_k^{j_k}$ is negative.

Recall that $r_k$ is defined as a linear combination of
$k$-standard monomials. Let ${\bf m}$ be a positive monomial of
$r_k$.
Even if ${\bf m} \in L$ is not a polynomial, the negative monomials of
${\bf m}$ are smaller than the largest negative monomial of
$g_k^{a_k}$ since $\deg_g({\bf m}) < \deg_g(g_k^{a_k})$ and
$\gap({\bf m}) \geq \gap(g_k)$.
So if e.g. $r_k$ does not contain negative
$k$-standard monomials then $\gap(g_{k+1}) = \gap(g_k)$.

In what follows $j$-standard monomials are ordered lexicographically by
their $g$-degree and $f$-degree, i.e. ${\bf m_i} < {\bf m_k}$ if $\overline{{\bf m_i}} < \overline{{\bf m_k}}$.
This order is well defined since $\overline{{\bf m}}$ determines ${\bf m}$ by Remark to Lemma 3.

To make reading less unpleasant we consider two cases: (i) $\gap(g_k) < \gap(g_{k-1})$ and (ii) $\gap(g_k) = \gap(g_{k-1})$.

\smallskip

(i) $\gap(g_k) < \gap(g_{k-1})$. Since $g_k = g_{k-1}^{a_{k-1}} -
r_{k-1}$ and $\gap(g_{k-1}^{a_{k-1}}) = \gap(g_{k-1}) > \gap(g_k)$
we can conclude that the largest negative monomial of $r_{k-1}$ is
larger than negative monomials of $g_{k-1}^{a_{k-1}}$. Since all
$k-1$-standard monomials have different $g$-degrees this monomial is
$\overline{\nu_{k-1}}$ for the largest negative $k-1$-standard
monomial $\nu_{k-1}$ of $r_{k-1}$. So $\gap(g_k) = \overline{g_{k-1}^{a_{k-1}}} \div
\overline{\nu_{k-1}}$.

Next, $g_{k+1} = (g_{k-1}^{a_{k-1}} - r_{k-1})^{a_k} - r_k =
g_{k-1}^{a_{k-1}a_k} - R_k - r_k$. Since $\deg_g(R_k) <
\deg_g(g_{k+1})$ we know that $R_k \in V_k$ (see Lemma 3).
Present $R_k$ through the standard basis as a sum of
$k$-standard monomials.

The largest negative $k$-standard monomial in $R_k$ turns out to be $\nu_{k-1}g_k^{a_k - 1}$.
Indeed
$\gap(g_{k-1}^{a_{k-1}a_k} - R_k) = \gap(g_k^{a_k}) = \gap(g_k) <
\gap(g_{k-1})$ and $\gap(g_{k-1}^{a_{k-1}a_k}) = \gap(g_{k-1})$; hence
the largest negative monomial of
$g_{k-1}^{a_{k-1}a_k}$ is smaller than the largest negative
monomial $\mu$ of $R_k$. Therefore $\overline{g_{k-1}^{a_{k-1}}} \div
\overline{\nu_{k-1}} =
\gap(g_k) = \overline{g_{k-1}^{a_{k-1}a_k}} \div \overline{\mu}$.
Since $\overline{g_{k-1}^{a_{k-1}}} = \overline{g_k}$ we
have $\overline{\mu} =
\overline{g_k^{a_k - 1}} \overline{\nu_{k-1}}$ and a $k$-standard monomial $\mu =
\nu_{k-1}g_k^{a_k - 1}$.

Let us compute its $z$-degree:
$\deg(\nu_{k-1}g_k^{a_k - 1}) = \deg(\nu_{k-1}) + (a_k - 1)m_k >
a_km_k$ because $\nu_{k-1}$ is a $k-1$-standard monomial of $r_{k-1}$
and
$\deg(\nu_{k-1}) > m_k = \deg(g_k)$. But $\deg(r_k) = a_km_k$ and
all $k$-standard monomials in $r_k$ have $z$-degree not exceeding
$a_km_k$. So $\nu_{k-1}g_k^{a_k - 1}$ is not a summand of $r_k$
and cannot be canceled.

\smallskip

(ii) $\gap(g_k) = \gap(g_{k-1})$. Since $\gap(g_0) = \infty$ and
$\gap(g_k) < \infty$ we can find such a $p$ that $\gap(g_k) =
\gap(g_{k-1}) = \dots = \gap(g_p) < \gap(g_{p-1})$. Just as above,
$g_{k+1} = g_{p-1}^{a_{p-1} \dots a_k} - R_k - r_k$ where $R_k \in
V_k$. Since $\gap(g_{p-1}^{a_{p-1} \dots a_k}) = \gap(g_{p-1}) >
\gap(g_{p-1}^{a_{p-1} \dots a_k} - R_k) = \gap(g_k) = \gap(g_p)$
we can conclude that the maximal negative $k$-standard monomial in
the standard representation of $R_k$ is $\nu_{p-1}g_p^{a_p -
1}\dots g_k^{a_k - 1}$, where $\nu_{p-1}$ is the largest negative
$p-1$-standard monomial in $r_{p-1}$. But $\deg(\nu_{p-1}g_p^{a_p
- 1}\dots g_{k-1}^{a_k - 1}) = \deg(\nu_{p-1}) + (a_p - 1)m_p +
\dots + (a_k - 1)m_k > a_km_k = \deg r_k$ since $\deg(\nu_{p-1}) >
m_p$ and $a_jm_j > m_{j+1}$. So again this monomial cannot be
canceled by a monomial from $r_k$.$\Box$\\

\textbf{Remark.} Negative powers of $f$ can appear in the finite characteristic case because though the function gap satisfies properties (a) and (c), property (b) should be modified. If $h$ is monic in $g, \ \ \chr(K) = p \neq 0$, and $d = p^\alpha d_1$ where $(p, d_1) = 1$ then $\gap(h^d) = (\gap(h))^{p^\alpha} \geq \gap(h)$.
$\Box$\\

If $\chr(K) = 0$ then, by Lemma 6, $B_s \subset K[f, g]$ and $h \in K[f, g]$ can be presented as a sum $h = \sum\limits_{{\bf m_j} \in B_s} \delta_{\bf j}{\bf m_j}$ where $\delta_{\bf j}(f) \in K[f]$.
(A similar description of $K[f,g]$ is obtained in [SU] when $f, \ g \in K[z_1, z_2, \dots, z_t]$ and are algebraically independent.)
Since the degrees of different $s$-standard monomials from $B_s$ are different $\mod\, n$ (see the proof of Lemma 3 (a)), the semigroup $\Pi(f, g)$ of degrees of non-zero elements of the subalgebra  $K[f,
g]$ is spanned by $n,\ m_0, ..., m_s$, i.e. $\Pi(f, g) = \Pi_s = \spn\{n,\ m_0, ..., m_s\}$.

If $1 \in \Pi(f, g)$ then the smallest of $n,\ m_0, ..., m_s$ is
$1$. If $m_i = 1$ then $d_i = 1$. As we observed, $d_{i+1} < d_i$, hence $i = s$ and $1 = m_s = d_s$.
Now we can prove by (reverse) induction that $d_j \in \Pi_j = \spn\{n, m_0, \dots, m_{j+1})$ for $j
\geq 0$. Assume that $d_{j+1} \in \Pi_{j+1}$ and $j > -1$. Since $d_{j+1} = (n, m_0, \dots, m_{j+1}) \in \Pi_{j+1}$ it is a linear combination of $\{n,\ m_0, ..., m_{j+1}\}$ with non-negative coefficients and $d_{j+1} = \min(n,\ m_0, ..., m_{j+1})$. If this minimum is $m_i$ where $i < j+1$ (here $n = m_{-1}$) then $d_i \leq m_i = d_{j+1}$ which is impossible because $d_{i+1} < d_i$ for $0 < i < s$. Therefore $m_{j+1} = d_{j+1}$ and ${d_j\over d_{j+1}}m_{j+1} =
d_j \in \Pi_j$. So $d_0 \in \Pi_0$ which proves the AMS. $\Box$

\bigskip

Also a beautiful result of David Richman that either ${n\over m_0}$ or ${m_0\over n}$ is an integer if $K[f, g]$ contains an element $h$ with the degree $d_0 = (n, m_0)$ (see [Ri], Proposition 1) follows from the presentation of $K[f, g]$ trough the standard monomials. Indeed, $d_0 = am_0 + bn$ where $0 \leq a < a_0$ and the standard monomial which has the degree $d_0$ must be $f^bg_0^a$. Since $b \geq 0$ either $n = d_0$ or $m_0 = d_0$.

\bigskip

\textbf{Remark.} If $\chr(K) = p$ and $d_0 = (n, m)$ is not divisible by $p$ the proof above is applicable verbatim: just assume that $m \not\equiv 0 {\pmod p}$ (switching $f$ and $g$ if necessary); then $a_i \not\equiv 0 {\pmod p}$ for $0 \leq i \leq s$, and all $g_i$ are polynomials of $f$ and $g$ since $\gap(g_i^{a_i}) = \gap(g_i)$. $\Box$\\

\textbf{Conclusion.}\\

In fact we proved a bit more: if $1 \in \Pi(f, g)$ then all ${m_i\over m_{i+1}}, \ i = 0, 1, \dots, s - 1$ are integers as well as ${n\over m_0}$ or ${m_0\over n}$.
We can call such a sequence \textit{1-admissible}. It is easy to show that any
$1$-admissible sequence can be realized by a pair of polynomials.

\bigskip

\textbf{Question.} Assume that $d$ is the smallest positive number in $\Pi(f, g)$. Describe all pairs $f, \ g$ for which this condition is satisfied.

\bigskip

If $d = 2$ and up to a change of variable $K[f, g] = K[z^2]$  then the question is already answered by the AMS Theorem. Another possibility is $f = zh(z^2), \ g = z^2$ where $\deg(h) > 1$. By the Richman's result mentioned above if $(n, m)$ is divisible by 2 then $\min(n, m) = (n, m)$. In a more interesting case when $\min(n, m) \neq (n, m)$ and hence $(n, m)$ is not divisible by 2 we may assume that $n$ is odd and
show with the approach used above that a \textit{2-admissible} sequence should be given by $n = (2b_t + 1)\cdot \dots \cdot (2b_0 + 1) \cdot (2b_{-1} + 1), \ m_0 = 2(2b_t + 1)\cdot \dots \cdot(2b_0 + 1), \ m_1 = 2(2b_t + 1)\cdot \dots \cdot(2b_1 + 1), \dots, m_t = 2(2b_t + 1), \ m_{t+1} = 2$ where $b_i$ are positive integers.
The smallest
non-trivial example $9, 6, 2$ of a $2$-admissible sequence
is realized by polynomials $f = z^9 + 6z^5 + 6z, \ g_0 = z^6
+ 4z^2$ since $g_1 = g_0^3 - f^2 + 8g_0 = -4z^2$. This pair is unique up to a change of variable (and
multiplying polynomials by constants to make them monic). Wen-Fong Ke showed using computer that the sequences
$(15, 6, 2)$, $(21, 6, 2)$, $(27, 6, 2)$, and $(15, 10, 2)$ cannot be realized.

\bigskip

\textbf{Conjecture.} If $2$ is the smallest positive number in $\Pi(f, g)$ and $n > m$ is odd, $m > 2$ is even then $n = 9, \ m = 6$.\\

\textbf{Acknowledgements.}\\

The author is grateful to the Max-Planck-Institut f\"{u}r Mathematik in Bonn, Germany where the
work on this project has been started (see [ML]). He was also supported by an NSA grant
H98230-09-1-0008, by an NSF grant DMS-0904713, a Fulbright fellowship awarded by
the United States–-Israel Educational Foundation, and a FAPESP grant 2011/52030-5 awarded by the State of S\~{a}o Paulo, Brazil.

\begin{center}
{\bf References Sited}
\end{center}

\noindent [AM] S. Abhyankar; T. Moh,  Embedding of the line in
the plane. J. Reine Angew. Math. 276 (1975), 148--166.

\noindent [AO] N. A$'$Campo; M. Oka, Geometry of plane curves via Tschrinhausen resolution tower. Osaka
J. Math. 333 (1996), 1003–-1033.

\noindent [AB] E. Artal-Bartolo, Une démonstration géométrique du théorème d'Abhyankar-Moh. J. Reine Angew. Math. 464 (1995), 97–-108.

\noindent [Es] A. van den Essen, Polynomial automorphisms and the
Jacobian Conjecture, Progress in Mathematics, 190. Burkh\"{a}user
Verlag, Basel, 2000.

\noindent [Gu] R. Gurjar,  A new proof of the Abhyankar-Moh-Suzuki theorem. Transform. Groups 7 (2002), no. 1, 61–-66.

\noindent [GM] R. Gurjar; M. Miyanishi,  On contractible curves in the complex affine plane. Tohoku Math.
J. 48(2) (1996), 459–-469.

\noindent [Ka] M. Kang, On Abhyankar-Moh's epimorphism theorem.
Amer. J. Math. 113 (1991), no. 3, 399–-421.

\noindent [ML] L. Makar-Limanov, A new proof of the Abhyankar-Moh-Suzuki theorem.
MPIM2005--77.

\noindent [Mi] M. Miyanishi,  Analytic irreducibility of certain curves on a nonsingular affine surface. In:
Proc. Int. Symp. in Algebraic Geometry, Kyoto 1977. Kinokuniya, Tokyo, 1978, pp. 575–-587

\noindent [PR] B. Peskin; D. Richman, A method to compute minimal polynomials.
SIAM J. Algebraic Discrete Methods 6 (1985), no. 2, 292–-299.

\noindent [Ri] D. Richman,  On the computation of minimal polynomials. J. Algebra
103 (1986), no. 1, 1–-17.

\noindent [Ru] L. Rudolph, Embeddings of the line in the plane. J. Reine Angew. Math. 337 (1982), 113–-118.

\noindent [SU] I. Shestakov; U. Umirbaev, Poisson brackets and two-generated subalgebras of rings of polynomials. J. Amer. Math. Soc. 17 (2004), no. 1, 181–-196.

\noindent [Su] M. Suzuki, Propi\'et\'es topologiques des
polyn\^omes de deux variables complexes, et automorphismes
alg\'earigue de l'espace $ C^2$. J. Math. Soc. Japan, 26 (1974),
241--257.

\noindent [Zo] H. \.{Z}o\l\c{a}dek, A new topological proof of the Abhyankar-Moh
theorem. Math. Z. 244 (2003), no. 4, 689–-695.

\vskip 4cm

\noindent Department of Mathematics \& Computer Science, the Weizmann Institute of Science,
Rehovot 76100, Israel;\\ Department of Mathematics, Wayne State University,
Detroit, MI 48202, USA; \\ Department of Mathematics, University of Michigan,
Ann Arbor, MI 48109, USA.
\vskip 0.5cm
\noindent \textit{E-mail address}: lml@math.wayne.edu\\

\end{document}